\documentclass[american]{amsart}
\usepackage[T1]{fontenc}
\usepackage[latin1]{inputenc}
\usepackage{babel}

\makeatletter

\theoremstyle{plain}    
\newtheorem{thm}{Theorem}[section]
\numberwithin{equation}{section} 
\numberwithin{figure}{section} 
\theoremstyle{plain}    
\newtheorem{lem}[thm]{Lemma} 
\theoremstyle{plain}    
\newtheorem{prop}[thm]{Proposition} 
\theoremstyle{plain}    
\theoremstyle{remark}
\newtheorem{rem}[thm]{Remark}
\theoremstyle{remark}    
\newtheorem{notation}[thm]{Notation} 
\newtheorem{assumptions}[thm]{Assumptions} 
\newtheorem{notation-assumptions}[thm]{Notation-Assumptions} 

\usepackage[arrow,matrix,curve,ps]{xy}
\CompileMatrices

\newcommand{\basefield}{k}

\renewcommand{\O}{\mbox{$\mathcal{O}$}}
\renewcommand{\P}{\mbox{$\mathbb{P}$}}

\newcommand{\C}{\mbox{$\mathbb{C}$}}

\newcommand{\sheafspec}{\mbox{\bf Spec}}
\newcommand{\mc}{\mathcal}

\DeclareMathOperator{\Aut}{Aut}
\DeclareMathOperator{\Chow}{Chow}
\DeclareMathOperator{\codim}{codim}

\DeclareMathOperator{\Hilb}{Hilb}
\DeclareMathOperator{\Hom}{Hom}
\DeclareMathOperator{\Image}{Image}
\DeclareMathOperator{\locus}{locus}
\DeclareMathOperator{\Pic}{Pic}

\DeclareMathOperator{\RatCurves}{RatCurves}
\DeclareMathOperator{\red}{red}
\DeclareMathOperator{\Reg}{Reg}
\DeclareMathOperator{\Sing}{Sing}
\DeclareMathOperator{\Singni}{Sing,ni}

\DeclareMathOperator{\Spec}{Spec}

\DeclareMathOperator{\Univ}{Univ}

\makeatother

\begin{document}

\title{Families of singular rational curves}

\date{\today }

\author{Stefan Kebekus}

\keywords{Rational curves of minimal degree, Fano Manifolds,
Bend-and-Break, Characterization of \protect $ \P _{n}\protect $}

\subjclass{Primary 14J45, Secondary 14C15}

\thanks{The author gratefully acknowledges support by a
Forschungsstipendium of the Deutsche Forschungsgemeinschaft}

\address{Stefan Kebekus, Institut für Mathematik, Universität
Bayreuth, 95440 Bayreuth, Germany}

\curraddr{Research Institute for Mathematical Sciences, Kyoto
University, Kyoto 606-8502, Japan}

\email{stefan.kebekus@uni-bayreuth.de}

\urladdr{http://btm8x5.mat.uni-bayreuth.de/\~{}kebekus}

\maketitle
\tableofcontents

\section{Introduction}

This work is concerned with the study of algebraic families of
rational curves. The main theorem yields a splitting-criterion for
families of singular curves. In some cases, this effectively
complements the bend-and-break argument which appears in the work of
Mori. We apply this result to projective varieties $X$ which are
covered by a family of rational curves of minimal degrees. More
precisely, we fix a general point $x\in X$ and prove in
theorem~\ref{thm:bounds_for_sing} that the subfamily of singular
rational curves which pass through $x$ is at most
one-dimensional. Furthermore, we describe the singularities of the
curves.

This has further consequences: first, we show that the tangent map
which sends a curve through $x$ to its tangent direction in
$\P(T_X^*|_x)$ is a finite morphism. We believe that this will be
useful in the further study of Fano manifolds with Picard number one.

For the second application, recall Mori's bend-and-break which asserts
that if $x,y\in X$ are two general points, then there are at most
finitely many curves in the family which contain both $x$ and $y$. In
this work we shed some light on the question as to whether two
sufficiently general points actually define a unique curve.

Finally, we give a characterization of the projective space which, in
characteristic 0, improves on some of the known generalizations of
Kobayashi-Ochiai's theorem.

Throughout this paper, unless otherwise noted, we work over an
algebraically closed field $\basefield$ of arbitrary characteristic
and use the standard notation of algebraic geometry. For spaces of
rational curves, our principal reference is \cite{K96}.

\subsubsection*{Acknowledgement }

The paper was written while the author enjoyed the hospitality of RIMS
in Kyoto.  The author is grateful to Y.~Miyaoka for the invitation and
to the members of that institute for creating a stimulating
atmosphere. He would like to thank S.~Helmke, J.-M.~Hwang, S.~Kovács,
Y.~Miyaoka and S.~Mori for a number of discussions on the subject.

\section{Families of singular curves}
\label{Sec:DefoSect}

The subject of the present section is a splitting criterion for
algebraic families of singular rational curves. Before stating the
main result, however, we need to introduce notation. Without this, the
formulation would be difficult.

\begin{notation}
Let $C$ be a curve and $\eta :\tilde{C}\to C$ the normalization. We
say that $C$ is ``immersed'' if $\eta $ has rank one at all points of
$\tilde{C}$.
\end{notation}
\begin{notation}
\label{Notation:genl_wrt_L}Let $C$ be an irreducible and reduced 
singular rational curve and $L\in \Pic (C)$ a line bundle of degree
$k>0$.  Let $(C_{i})_{i\in \{0,1\}}$ be the nodal and cuspidal plane
cubic. We say that a smooth point $x\in C$ is ``general with respect
to $L$'', if $\beta ^{*}_{i}(L)\not \cong \O _{C_{i}}\left( k\beta
^{-1}_{i}(x)\right) $ for all birational morphisms $\beta
_{i}:C_{i}\to C$ and for all $i\in \{0,1\}$.
\end{notation}

\begin{lem}
\label{lem:genl_pt_is_L_genl}
If $C$ and $L$ are as above, then the set of points which are general
with respect to $L$ is Zariski-open in $C$.
\end{lem}

\begin{proof}
As a first step, note that if $x\in C$ is a smooth point, $g\in \Aut
(C_i)$ an automorphism and $\beta_i:C_i \to C$ a birational morphism,
then it is clear that
$$
(g\circ \beta_i)^{-1}(L)\cong 
\O_{C_i}\left( k(g\circ \beta_i)^{-1}(x)\right) 
$$
if and only if $\beta_i^{-1}(L)\cong \O _{C_i}\left( k\beta_i^{-1}(x)\right)$.

To conclude, it suffices to note that ---up to composition with
automorphisms--- there are only finitely many birational morphisms
$C_i \to C$. Recall that for a given bundle $L_i\in \Pic (C_i)$, we
have that $\O_{C_i}(ky) \not \cong L_i$ for all but finitely many
smooth points $y\in C_i$.
\end{proof}

The following is the main result of this section and the technical
core of this paper.

\begin{thm}
\label{Thm:Singular_b_and_b}
Let $\pi :X\to Y$ be a projective morphism between proper
positive-dimensional varieties such that the general fiber is an
irreducible and singular rational curve. Assume that $\pi$ is
equi-dimensional and that there exists a section $\sigma _{\infty
}\subset X$ and a morphism $\gamma :X\to Z$ which maps $\sigma
_{\infty }$ to a point and is finite on the complement of $\sigma
_{\infty }$. Let $L\in \Pic (X)$ be a line bundle whose restriction to
$\pi $-fibers is of positive degree.  Assume furthermore that one of
the following holds:
\begin{enumerate}
\item all $\pi $-fibers are immersed curves 
\item no $\pi $-fiber is immersed and there exists a closed point
$y\in Y$ with reduced and irreducible fiber $X_{y}:=\pi ^{-1}(y)$ such
that $\sigma _{\infty }\cap X_{y}$ is a smooth point which is general
with respect to $L|_{X_{y}}$.
\item the restriction of $L$ to $\pi $-fibers is of degree 2 and there
exists a closed point $y\in Y$ with reduced and irreducible fiber
$X_{y}:=\pi ^{-1}(y)$ such that $\sigma _{\infty }\cap X_{y}$ is
general with respect to $L|_{X_{y}}$.
\end{enumerate}
Then there exists a point $y\in Y$ such that $\pi ^{-1}(y)$ is not
irreducible or not generically reduced.
\end{thm}

The remaining part of section~\ref{Sec:DefoSect} is devoted to the
proof of theorem~\ref{Thm:Singular_b_and_b}. The reader who is
predominantly interested in the applications may want to skip the rest
of this section on first reading. We remark that statement (1) has
already been considered by several authors, e.g.~\cite[sect.~2]{CS95}.

In order to prove the theorem, we argue by contradiction. More
precisely, we assume the following throughout the present
section~\ref{Sec:DefoSect}.

\begin{assumptions}
\label{Ass:Main_thm_wrong}
Assume that $\pi :X\to Y$ is a family as described in
theorem~\ref{Thm:Singular_b_and_b}, but all fibers of $\pi $ are
irreducible and generically reduced. In particular, the reduction of
any $\pi $-fiber is a singular rational curve. Since the statement of
theorem~\ref{Thm:Singular_b_and_b} is stable under finite base change
and under restriction to positive-dimensional subvarieties of $Y$, we
assume without loss of generality that $Y$ is a smooth curve.
\end{assumptions}

\subsection{Reduction to families of plane cubics }

As a first step in the proof of theorem~\ref{Thm:Singular_b_and_b}, we
show that the family $X\to Y$ can be replaced by a family where every
fiber is isomorphic to a plane cubic. Although we need this only over
smooth curves, we formulate more generally for arbitrary normal
bases. The succeeding lemma is the key to the reduction.

\begin{lem}\label{Lem:Cover_by_plane_cubics}
Consider the following diagram of surjective projective morphisms:
$$
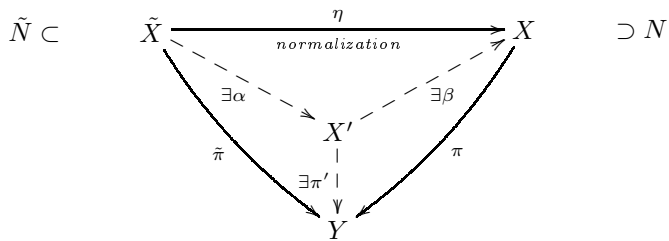
\xymatrix{{\tilde N} \subset & {\tilde X}
\ar@{-->}[rrd]_{\exists \alpha} \ar[rrrr]^{\eta}_{\txt{\tiny
normalization}} \ar@/_/ [rrdd]_ {\tilde \pi} & & & & {X}
\ar@/^/[lldd]^{\pi} & {\supset N} \\ & & & {X'}
\ar@{-->}[rru]_{\exists \beta} \ar@{-->}[d]_{\exists \pi'} & & \\& & &
{Y}& &}
$$ 
where $\pi $ is projective, $Y$ is a smooth curve and the following holds
\begin{itemize}
\item $\tilde \pi$ is a $\P_1$-bundle, i.e.~a smooth morphism whose fibers
are isomorphic to $\P_1$

\item all $\pi $-fibers are irreducible and generically reduced
singular rational curves 

\item there exists a section $N\cong Y$ and a subscheme
$\tilde{N}\subset \eta ^{-1}(N)$ so that for all closed points $y\in
Y$, the scheme-theoretic intersection $\tilde{\pi }^{-1}(y)\cap N$ is
a zero-dimensional subscheme of length 2
\end{itemize}
Then there exists a factorization $\eta =\beta \circ \alpha $ via a
variety $X'$ such that all fibers of $\pi '$ are rational curves with
a single cusp or node.
\end{lem}
\begin{proof}
Set
$$
\mc A:=\Image \left( \eta _{*}(\mc J_{\tilde{N}})\oplus \O_X\to
\eta _{*}(\O _{\tilde{X}})\right) \subset \eta _{*}(\O _{\tilde{X}}),
$$
where $\mc J_{\tilde{N}}\subset \O _{\tilde{X}}$ is the ideal sheaf of
$\tilde{N}$. It follows immediately that $\mc A$ is a coherent sheaf
of $\O _{X}$-algebras. Define $X':=\sheafspec (\mc A)$. The existence
of $\alpha $, $\beta $ and $\pi '$ follows by construction.

In order to see that fibers of $\pi '$ are of the desired type, let
$y\in Y$ be an arbitrary closed point. We are finished if we show that
the fiber $(\pi ')^{-1}(y)$ has a single singularity which is a simple
node or cusp.

After replacing the base $Y$ by an affine neighborhood of $y$ and
performing a base change, if necessary, we may assume that there
exists a relatively ample divisor $D\subset X$ which intersects every
$\pi $-fiber in a single smooth point. We may furthermore assume that
$\tilde{X}\cong Y\times \P _{1}$.  Write $U:=X\setminus D$, write
$\tilde{U} := \eta^{-1}(U)$ and note that both $U$ and $\tilde{U}$ are
affine. By choosing a bundle coordinate $z$ on $\tilde{X}$, we may
write $U\cong \Spec (R)$, $Y\cong \Spec (S)$ and $\tilde{U}\cong \Spec
(S\otimes_\basefield \basefield[z])$.

Because $N\cong H$ we can decompose $R=\mc J_{N}(U)\oplus \pi ^{\#}(S)$.
Since $\eta ^{\#}\left( \mc J_{N}(U)\right) \subset \mc J_{\tilde{N}}$
by construction, we have the equation of rings
$$
\begin{array}{ccc}
\mc A(U) & = & \Image \left( \mc J_{\tilde{N}}(\tilde{U})\oplus R\to
 S\otimes \basefield[z]\right) \\ 
& = & \Image \left( \mc J_{\tilde{N}}(\tilde{U})\oplus S\to S\otimes \basefield[z]\right)
\end{array}
$$
If $z$ is chosen properly, tensoring with $\basefield(y)$ yields
$$
\begin{array}{ccc}
\mc A(U)\otimes \basefield(y) & = & \Image \left( \mc
J_{\tilde{N}}(\tilde{U})\otimes \basefield(y)\oplus \basefield(y)\to
\basefield(y)[z]\right) \\ 
& = & (z^{2}-c)+\basefield(y)\subset \basefield(y)[z]
\end{array}
$$
for a number $c\in \basefield(y)$. An elementary calculation shows that this
ring is generated by the elements $z^{2}-c$, $z(z^{2}-c)$ and the
constants $\basefield(y)$. Therefore
$$
\begin{array}{ccc}
\beta ^{-1}\left( \pi ^{-1}(y)\cap U\right) & = & \Spec \basefield(y)\left[
 z^{2}-c,z(z^{2}-c)\right] \\ 
& = & \Spec \basefield(y)[a,b]\, /\, \left(
 b^{2}-a^{2}(a+c)\right)
\end{array}
$$
which defines a cusp if $c=0$ and a node otherwise.
\end{proof}

In the proposition above, the assumption that $\tilde \pi$ is a
$\P_1$-bundle is automatically satisfied in characteristic zero. In
positive characteristic, it can be necessary to perform a finite base
change before the normalization yields a bundle. See \cite[II.2]{K96}
for a detailed account of this phenomenon.

\begin{prop}
\label{Prop:Cover_by_plane_cubics}
If $\pi :X\to Y$ satisfies the assumptions~\ref{Ass:Main_thm_wrong},
then there is a family $\pi ':X'\to Y'$ which satisfies
assumptions~\ref{Ass:Main_thm_wrong} as well, and all fibers are
rational curves with a single node or cusp.
\end{prop}

\begin{proof}
The strategy of this proof is to find a sequence of base changes which
modify $X$ and $Y$ so that lemma~\ref{Lem:Cover_by_plane_cubics} can
be applied.

First, after finite base change, we may assume that the normalization
$\tilde X$ is a $\P_1$-bundle over $Y$. Recall that smooth morphisms
are stable under base change. Thus, even after further base changes,
the normalization of the pull-back of $X$ will still be a
$\P_1$-bundle over the base.

Now let $X_{\Sing,\pi}\subset X$ be the singular locus of $\pi$-fibers
and let $N\subset X_{\Sing,\pi}$ be an irreducible component which
maps surjectively onto $Y$. Since $Y$ is assumed to be a smooth curve,
$N$ is actually finite over $Y$. If $\pi|_N$ is not isomorphic,
perform a base change. Thus, we assume that $N\cong Y$.

As a next step, let $\eta :\tilde{X}\to X$ be the normalization and
consider the relative Hilbert-scheme $\pi_H:\Hilb_2(\eta
^{-1}(N)/Y)\to Y$ of zero-dimensional subschemes of length 2 in $\eta
^{-1}(N)$ over $Y$.  Recall the fact that taking the Hilbert-scheme
commutes with base change (see e.g. \cite[I.1.4.1.5]{K96}). In our
setup this means that if $Y'\to Y$ is a morphism, then
$$
\Hilb_2(\eta ^{-1}(N)\times_Y Y'/Y')\cong \Hilb_2(\eta ^{-1}(N)/Y)\times_Y Y'.
$$
This has two important consequences.

First, by choice of $N$, if $y\in Y$ is a general point, then $\eta
^{-1}(N)\cap \tilde{\pi }^{-1}(y)$ is zero-dimensional of length at
least two, so that $\pi_H^{-1}(y)$ is not empty. It follows that
$\pi_H$ is surjective.

Second, let $V\subset (\Hilb_2 (\eta ^{-1}(N)/Y))_{\red }$ be a
subvariety which is finite over $Y$. Let $Y'$ be the normalization of
$V$ and perform another base change. Since taking $\Hilb$ commutes
with base change, $Y'$ can be seen as a subscheme of $\Hilb_2(\eta
^{-1}(N)\times_Y Y'/Y')$ and therefore defines a subscheme
$\tilde{N}\subset \widetilde{X\times_Y Y'}=\tilde{X}\times_Y Y'$.
Thus, all the prerequisites of lemma~\ref{Lem:Cover_by_plane_cubics}
are fulfilled, and we can apply that lemma in order to obtain
$X'$. Because our construction involves only finite base change, it is
clear that $X'$ satisfies the requirements of
theorem~\ref{Thm:Singular_b_and_b} if and only if $X$ does.
\end{proof}
\begin{rem}
\label{Rem:Pure_type_of_singularities}
Under the assumptions~\ref{Ass:Main_thm_wrong}, let $X_{y}$ be a
general $\pi $-fiber. If $X_{y}$ is not immersed, then it is clear
from the construction that one can choose $\pi ':X'\to Y'$ so that any
$\pi '$-fiber is isomorphic to a cuspidal plane cubic. Analogously, if
all $\pi $-fibers are immersed, the construction automatically yields
a family of curves where each fiber is isomorphic to a nodal plane
cubic.
\end{rem}

\subsection{Proof of theorem~\ref{Thm:Singular_b_and_b}}

We will now employ the line bundle $L$ in order to find two disjoint
sections $\tilde{\sigma }_0,\tilde{\sigma }_1 \subset X$ which cannot
be contracted.  The following elementary lemma says that this is not
possible.

\begin{lem}
\label{Lem:Sections_in_ruled_surfaces}
Let $\pi :\tilde{X}\to Y$ be a smooth minimal ruled surface. Assume
that there are three distinct sections $\tilde{\sigma }_{0}$,
$\tilde{\sigma }_{1}$ and $\tilde{\sigma }_{\infty }$ where
$\tilde{\sigma }^{2}_{\infty }<0$.  Then $\tilde{\sigma }_{0}$ and
$\tilde{\sigma }_{1}$ are not disjoint.
\end{lem}
\begin{proof}
Write $\tilde{\sigma }_{0}\equiv \tilde{\sigma }_{\infty }+a_{0}C$,
$\tilde{\sigma }_{1}\equiv \tilde{\sigma }_{\infty }+a_{1}C$ where
$\equiv $ denotes numerical equivalence and $C$ is a general fiber of
$\pi $. Since $\tilde{\sigma }_{0}$ and $\tilde{\sigma }_{1}$ are
effective, $\tilde{\sigma }_{0,1}.\tilde{\sigma }_{\infty }\geq 0$ and
thus $a_{0},a_{1}\geq -\tilde{\sigma }^{2}_{\infty }>0$. Therefore
$\tilde{\sigma }_{0}.\tilde{\sigma }_{1}=\tilde{\sigma }^{2}_{\infty
}+a_{0}+a_{1}>0$.
\end{proof}
With this preparation we can now finish the proof of
theorem~\ref{Thm:Singular_b_and_b}.  We stick to the
assumptions~\ref{Ass:Main_thm_wrong} and let $\pi :X\to Y$ be the
family of curves with a single cusp or node whose existence is
guaranteed by proposition~\ref{Prop:Cover_by_plane_cubics}. Let $\eta
:\tilde{X}\to X$ be the normalization and recall that the natural map
$\tilde{\pi }:\tilde{X}\to Y$ gives $\tilde{X}$ the structure of a $\P
_{1}$-bundle.

\subsubsection{Proof of case (1) of theorem~\ref{Thm:Singular_b_and_b}}

In this setting, we may assume that all $\pi $-fibers are nodal plane
cubics.  After base change, we may assume that the preimage $\eta
^{-1}(X_{\Sing })$ contains two disjoint sections $\tilde{\sigma
}_{0}$ and $\tilde{\sigma }_{1}$.  If $\sigma _{\infty }\subset
X_{\Sing }$, then we can choose $\tilde{\sigma }_{0},\tilde{\sigma
}_{1}\subset \eta ^{-1}(\sigma _{\infty })$ and both $\tilde{\sigma
}_{0}$ and $\tilde{\sigma }_{1}$ can be contracted.  This is clearly
impossible. On the other hand, if $\sigma _{\infty }\not \subset
X_{\Sing }$, then $\tilde{X}$ contains the section $\eta ^{-1}(\sigma
_{\infty })$, which can be contracted, and the sections $\tilde{\sigma
}_{0},\tilde{\sigma }_{1}\subset \eta ^{-1}(X_{\Sing })$ which are
disjoint. Since $\tilde{\sigma }_{0}$ and $\tilde{\sigma }_{1}$ are
both different from $\eta ^{-1}(\sigma _{\infty })$, this contradicts
lemma~\ref{Lem:Sections_in_ruled_surfaces}.$\hfill \qed $

\subsubsection{Proof of case (2) of theorem~\ref{Thm:Singular_b_and_b} }

Let $y\in Y$ be any closed point and set $X_{y}:=\pi ^{-1}(y)$ and
$k:=\deg L|_{X_{y}}$. By remark~\ref{Rem:Pure_type_of_singularities},
we may assume that $X_{y}$ is a rational curve with a single cusp. But
then there exists a unique smooth point $x_{y}\in X_{y}$ such that $\O
_{X_{y}}(kx_{y})\cong L|_{X_{y}}$.  A calculation of the deformation
space of the cuspidal plane cubic yields that $\pi $ is a locally
trivial fiber bundle. In particular, by taking the union of the
$x_{y}$ we obtain a section $\sigma _{0}\subset X_{\Reg }$ which does
not meet the singular locus of $X$. Set $\tilde{\sigma }_{0}:=\eta
^{-1}(\sigma_0)$, set $\tilde{\sigma}_1:=\eta ^{-1}(X_{\Sing })$
and note that $\tilde{\sigma }_{0}$ and $\tilde{\sigma }_{1}$ are
disjoint. By the assumption that $\sigma _{\infty }\cap X_{y}$ is a
smooth point which is general with respect to $L|_{X_{y}}$, it follows
that $\tilde{\sigma }_{0}$, $\tilde{\sigma }_{1}$ and $\eta
^{-1}(\sigma _{\infty })$ are distinct. Again, this contradicts
lemma~\ref{Lem:Sections_in_ruled_surfaces}.$\hfill \qed $

\subsubsection{Proof of case (3) of theorem~\ref{Thm:Singular_b_and_b}}

We may assume without loss of generality that there exist points $y\in
Y$ such that $\pi ^{-1}(y)$ is a nodal curve. Otherwise, we are in
case (2) of the theorem.

Let $y\in Y$ be any closed point and set $X_{y}:=\pi ^{-1}(y)$. Since
$X_{y}$ is a rational curve with a single node or cusp, and
$L|_{X_{y}}$ is a line bundle of degree two, $L|_{X_{y}}$ is
basepoint-free and induces a 2:1 cover
$$
\Gamma _{y}:X_{y}\to \P(H^{0}(L|_{X_{y}})^{*})\cong \P _{1}.
$$
If $X_{y}$ is a nodal curve, and $\eta :\P _{1}\to X_{y}$ the
normalization, we know that the branch locus of $\Gamma _{y}\circ \eta
$ consists of two distinct points $x_{0},x_{1}$ which are not
contained in the preimage of the singularity: $x_{i}\not \in \eta
^{-1}((X_{y})_{\Sing })$. We extend $\Gamma _{y}$ to a global map
$\Gamma $: 
$$
\xymatrix{{\tilde
X} \ar[rr]^{\eta}_{\txt{\tiny normalization}} \ar@/_/[rrd]_{\tilde
\pi} & & {X} \ar[r]^-{\Gamma} \ar[d]^{\pi}& {\P(\pi_*(L)^*)}
\ar@/^/[ld] \\ & & {Y} & }
$$ 
Note that $\pi _{*}(L)$ is locally free of rank 2. Recall that
$\tilde{X}$ is a minimal smooth ruled surface over $Y$ and let
$D\subset \tilde{X}$ be the (reduced) branch locus of $\Gamma \circ
\eta $. Then $D$ intersects every $\tilde{\pi }$-fiber in exactly 2
distinct points.

After performing another base change, if necessary, we may assume
without loss of generality that $D$ is reducible. Write
$D=\tilde{\sigma }_{0}\cup \tilde{\sigma }_{1}$ and note that
$\tilde{\sigma }_{0}\cap \tilde{\sigma }_{1}=\emptyset $, i.e. that
$\tilde{\sigma }_{0}$ and $\tilde{\sigma }_{1}$ are disjoint
sections. But $\tilde{X}$ contains also a contractible section
$\tilde{\sigma }_{\infty }$, which is a component of the preimage of
$x\in X$. By assumption (3) of theorem~\ref{Thm:Singular_b_and_b},
these three sections are distinct, contradiction. This ends the proof
of theorem~\ref{Thm:Singular_b_and_b}.$\hfill \qed $

\section{Applications}

\subsection{Families of singular curves on projective varieties}

The purpose of this section is to give an estimate  of the dimension
of the space of singular curves through a general point and the
describe the singularities.

For the reader's convenience, we recall some facts about parameter
spaces for rational curves. See~\cite[II.2]{K96} for a detailed
account of these matters. Let $X$ be a projective variety and
$\Chow(X)$ its Chow-variety with universal family
$$
\Univ ^{\Chow }(X)\subset X\times \Chow (X).
$$
It can be shown that $\Chow (X)$ contains a quasi-projective
subvariety $V\subset \Chow (X)$ parameterizing irreducible and
generically reduced rational curves. 

In characteristic~0, define $\RatCurves^n(X)$ to be the normalization
of $V$, and note that the normalization of the universal family over
$\RatCurves^n (X)$ becomes a $\P_1$-bundle. In order to achieve the
same in positive characteristic, let $\RatCurves^n(X)$ be the
normalization of a canonically given finite cover of $V$.

Throughout the rest of this work we consider the following setup.

\begin{notation-assumptions}
\label{Ass:Properness}
Let $X$ is a projective variety (not necessarily normal) over
$\basefield$ and $H\subset \RatCurves ^{n}(X)$ an irreducible family
(not necessarily proper) of rational curves with universal family
$U\subset X\times H$ as follows
$$
\xymatrix{
{U} \ar[r]^{\iota} \ar[d]_{\pi} & {X} \\ 
H
}
$$ 
Assume that $\iota$ is dominant.  If $x\in X$ is any point, then 
define $H_{x}\subset H$ to be the subfamily of curves through $x$:
$H_{X}:=\pi (\iota ^{-1}(x))_{\red }$.  Write $U_x$ for the
restriction of the universal family and $\iota_x$, $\pi_x$ for the
restrictions of the canonical morphisms.

We assume throughout the present section that $H_x$ is proper for
general choice of $x\in X$. We define $\locus (H_x):=\iota (U_x)$ with
its reduced structure.
\end{notation-assumptions}

\begin{rem}
It is a main result of Mori theory that a family $H$ (``minimal
rational curves'') satisfying the assumptions always exists if $X$ is
a Fano manifold.
\end{rem}

\begin{thm}
\label{thm:bounds_for_sing}
Let $H^{\Sing} \subset H$ be the subfamily parametrizing singular
curves and let $L\in \Pic(X)$ be a line bundle whose restriction to
the curves is of positive degree. If $x\in X$ is a general point, then
the following holds
\begin{enumerate}
\item The subfamily $H^{\Sing}_x$ of singular curves through $x$ has
dimension at most one. The subfamily $H^{\Sing,x}_x \subset
H^{\Sing}_x$ of curves which are singular at $x$ is at most finite. If
$H^{\Sing,x}_x$ is not empty, then the associated curves are immersed.

\item If $L$ intersects the curves with multiplicity two, then
$H^{\Sing}_x$ is at most finite and $H^{\Sing,x}_x$ is empty.
\end{enumerate}
\end{thm}

For the applications it is important to keep in mind that the notion
of a general point depends on the choice of the line bundle $L$.

\begin{proof}

As a first step, we need to find an estimate for the dimension of the
subfamily of non-immersed curves. We let $H^{\Singni}\subset H^{\Sing
}$ be the closed subfamily of non-immersed curves and claim that
\begin{equation}
\label{eq:bound_for_singni}
\dim H^{\Singni}_{x} < 1.
\end{equation}
Indeed, if $\dim H^{\Singni}_{x}\geq 1$, then let $C\subset X$ be a
singular curve which corresponds to a general point of a component of
$H^{\Singni}$ which is of maximal dimension. By
lemma~\ref{lem:genl_pt_is_L_genl}, we find a smooth point $y\in C$
which is general with respect to $L|_C$. We remark that $C$ is not an
isolated point point of $H^{\Singni}_y$ and conclude by
theorem~\ref{Thm:Singular_b_and_b}.(2) that $H^{\Singni}_y$ cannot be
proper. But then $H_y$ cannot be proper, contrary to our
assumption. This shows the inequality~(\ref{eq:bound_for_singni}).

Next, we show that $\dim H^{\Sing}_x \leq 1$. We will argue by
contradiction and assume that $\dim H^{\Sing}_x \geq 2$. Now, if
$H'\subset H^{\Sing}_x$ is any curve which parameterizes nodal curves
in $X$, we can apply theorem~\ref{Thm:Singular_b_and_b}.(1) to the
family $H'$ to see that $H'$ is not proper. It follows that the
closure $\overline{H'}$ intersects $H^{\Singni}_x$ and therefore
$$
\codim_{H^{\Sing }_{x}}H^{\Singni}_{x}\leq 1.
$$
Hence our claim follows from the inequality~(\ref{eq:bound_for_singni}).

In order to show that $H^{\Sing,x}_x$ is at most finite, perform a
dimension count. It is clear that 
$$
\dim U^{\Sing} = \dim H^{\Sing} +1 \geq \dim H^{\Sing,x}_x + \dim X + 1,
$$
which in turn implies that general fibers $F$ of the natural
projection $U^{\Sing} \to X$ are of dimension $\dim F \geq
\dim H^{\Sing,x}_x + 1$. Now it suffices to note that the natural map $F\to
H^{\Sing}_x$ is finite, i.e. that $\dim F = \dim H^{\Sing}_{x}$ to
obtain that 
$$ 
1 \geq \dim H^{\Sing}_x \geq \dim H^{\Sing,x}_x + 1
$$ 
which yields the finiteness result.

With the inequality~(\ref{eq:bound_for_singni}), the same dimension
count, using the family $H^{\Singni}$, immediately shows that
$$ 
1 > \dim H^{\Singni}_x \geq \dim H^{\Singni,x}_x + 1 
$$
which means that all curves associated with $H^{\Sing,x}_x$ are
immersed. This ends the proof of statement (1).

To prove statement (2), we argue as above. Let $C\subset X$ be a
singular curve corresponding to a general point of $H^{\Sing}$ and
$y\in C$ general with respect to $L|_C$. By
theorem~\ref{Thm:Singular_b_and_b}.(3), the family $H^{\Sing}_y$
cannot be proper, contradiction.

\end{proof}

\subsection{Existence of a finite tangent morphism}

If $\dim H_x > 1$, then theorem~\ref{thm:bounds_for_sing}.(1) asserts
that the general curve associated with $H_x$ is smooth at $x$. If
$\tilde H_x$ is the normalization of $H_x$, we may therefore define
the ``tangent map''
$$
\tau_x : \tilde H_x \dasharrow \P(T^*_X|_x)
$$
which associates a curve through $x$ which is smooth at $x$ with its
tangent direction at $x$. Here we show that $\tau_x$ is always a
morphism.

The image of the tangent map has been studied extensively in a series
of papers by J.-M.~Hwang and N.~Mok. The authors informed us that
theorem~\ref{thm:finiteness_of_tau} can be used to give a different,
simpler proof of the deformation rigidity of hermitian symmetric
spaces.

\begin{thm}
\label{thm:finiteness_of_tau}
Under the assumptions~\ref{Ass:Properness}, if $x\in X$ is a general
point, then the tangent map $\tau_x$ is a finite morphism.
\end{thm}

\begin{proof}
If $f:\P_1 \to X$ is a morphism which is birational onto its image, if
$f(0)=x$ and the image $\Image(f)$ is a curve which is associated with
$\tilde H_x$, then it follows from theorem~\ref{thm:bounds_for_sing}.(1)
that $f$ is smooth in a neighborhood of $0\in \P_1$. Now, to conclude
that $\tau_x$ is a morphism, it suffices to realize that $\tau_x$ can
be written as a composition
$$
\xymatrix{
{\tilde H_x} \ar@/_/[rrd]_{\tau_x} \ar[rr]^(.3){\alpha} & & 
\Hom_{bir}(\P_1,X,0 \mapsto x)/\Aut(\P_1,0) \ar[d]_{\beta}\\
 & & {\P(T^*_X|_x)}
}
$$ 
where $\Hom_{bir}(\P_1,X,0\mapsto x)$ is the parameter space of
morphisms, and $\beta$ sends a morphism $f$ to the image of its
tangent map: $\beta:f \mapsto Tf(T_{\P_1}|_0)$. The existence of
$\alpha$ and the existence of the quotient space follows from 
universal properties ---see \cite[II.3]{K96}.

To show that $\tau_x$ is finite, we argue by contradiction. If
$\tau_x$ was {\em not} finite, we could find a curve $H'\subset \tilde
H_x$ such that $\tau_x(H')$ is a single point $v \in
\P(T^*_X|_x)$. After changing base, we can assume that $H'$ is smooth,
and consider the diagram:
$$
\xymatrix{
U' \ar[d]_{\pi'} \ar[r]^{\iota'} & X \\
H' 
}
$$ 
where $U'$ is the normalization of the universal family and hence a
$\P_1$-bundle over $H'$. There exists a section $\sigma_\infty \subset
U'$ which contracts to the point $x$, and the restriction of the
tangent morphism $T\iota'$ yields a morphism:
$$
T\iota' : N_{\sigma_\infty,U'} \to v\cong \basefield
$$ 
where $N_{\sigma_\infty,U'}$ is the normal bundle of $\sigma_\infty$
and $v$ a line in $T_X|_x$. But since the normal bundle is not
trivial, this map has to have a zero! Thus, there exist curves
associated with $H'$ which have non-immersed singularities at
$x$. This contradicts theorem~\ref{thm:bounds_for_sing}.(1), and we
are done.
\end{proof}

\subsection{Uniqueness of rational curves through 2 points}
\label{Sec:biratl_criterion}

In the setup given in~\ref{Ass:Properness}, the author conjectures
that the map $\iota_x$ is always generically one-to-one onto its
image. It is a direct consequence of
theorem~\ref{thm:bounds_for_sing}.(2) that this holds in the presence
of a line bundle of low degree.

\begin{thm}\label{Thm:families_of_curves}
If there exists a line bundle $L\in \Pic (X)$ intersecting the curves
with multiplicity 2, then $\iota_x$ is generically one-to-one onto its
image. In particular, if $y\in \locus (H_x)$ is a general point, then
there exists a unique curve in $H$ containing both $x$ and $y$.
\end{thm}
\begin{proof}
Following an argument of Miyaoka, $\iota_x$ is generically one-to-one
if $H^{\Sing,x}_x = \emptyset$, and the latter follows from
theorem~\ref{thm:bounds_for_sing}.(2). See \cite[V.3.7.5]{K96} for a
proof of Miyaoka's result.
\end{proof}

Remark that if $k$ is a field of characteristic zero, then
``generically one-to-one'' implies ``birational''. In finite
characteristic this is of course not necessarily so, and $\iota_x$
will certainly not be birational in general.

\subsection{Characterization of $\P_n$}
\label{Sec:Kobayashi-Ochiai}

In \cite{KO73}, Kobayashi and Ochiai characterized the projective
space as the only Fano-manifold $X$ whose canonical bundle is
divisible by $\dim X+1$. This result was generalized by several
authors, e.g.~in \cite{KS99}. We employ
theorem~\ref{Thm:families_of_curves} to give a characterization of the
projective space in terms of families of rational curves which
improves on the known results.

\begin{thm}
\label{Thm:Characerization}
Let $X$ be a normal projective variety defined over $\C$, $L\in \Pic
(X)$ and $H\subset \RatCurves^n(X)$ an irreducible component. Let
$x\in \locus (H)$ be a general closed point. Assume that $H_x$ is
proper and that $\locus (H_x)=X$. If $L.C=2$ for a curve $C\in H$,
then $X\cong \P_n$.
\end{thm}
By the classic argumentation of Mori, when $X$ is a Fano-manifold and
$-K_{X}.C>\dim X$ for all rational curves $C\subset X$, then a family
$H$ exists where $H_{x}$ is proper and $\locus (H_x)=X$.

\begin{proof}
Let $\tilde{U}_x$ be the normalization of the universal family
$U_x$ and consider the diagram
$$
\xymatrix{
{\tilde U_x} \ar[r]^{\tilde \iota_x} 
\ar[d]_{\txt{$\tilde \pi_x$\\\scriptsize $\P_1$-bundle}} & X \\ 
{\tilde H_x} \ar[r]^(.4){\tau_x} & {\P(T^*_X|_x)}
}
$$ 
It follows directly from theorem~\ref{Thm:families_of_curves} that
$\tilde \iota_x$ is birational and it follows immediately that
$\tau_x$ is birational as well. But then $\tau_x$ is a finite
birational morphism between normal spaces, and therefore
isomorphic. In particular, $\tilde{H}_x$ and $\tilde{U}_x$ are smooth.

The induced map $\tilde{\iota }_{x}:\tilde{U}_{x}\to X$ is an
isomorphism away from the section $D:=\tilde{\iota}^{-1}_x (x)$, which
is contracted: otherwise, Zariski's main theorem asserts that there
exists a point $x'\in X$ and a positive dimensional subfamily of
curves passing through both $x$ and $x'$. But Mori's bend-and-break
argument says that this cannot happen if $H_x$ is unsplit. In
particular, since $\tilde{U}_x$ is smooth it follows that $X$ is
smooth.

In this setting, an elementary theorem of Mori yields the claim. See
\cite[V.3.7.8]{K96}.
\end{proof}

\end{document}